\documentclass{article}

\newtheorem{fed}{\textbf{Definition}}[section]

\newtheorem{lemma}[fed]{\textbf{Lemma}}

\newtheorem{prop}[fed]{\textbf{Proposition}}
\newtheorem{cor}[fed]{\textbf{Corollary}}
\usepackage{amssymb,bbm,graphicx,epsfig,psfrag,epic,eepic,latexsym}
\usepackage{amsmath}
\usepackage{mathrsfs}
\usepackage{color}
\usepackage{authblk}

\title{On doubly symmetric periodic orbits}
\author{Urs Frauenfelder, Agustin Moreno} 

\newcommand{\Addresses}{{
  \bigskip
  \footnotesize

  U.~Frauenfelder, \textsc{Augsburg Universität, Augsburg, Germany}\par\nopagebreak
  \textit{E-mail address}: \texttt{urs.frauenfelder@math.uni-augsburg.de}

  \medskip

  A.~Moreno, \textsc{Institute for Advanced Study, Princeton NJ, USA/ Heidelberg Universität, Heidelberg, Germany}\par\nopagebreak
  \textit{E-mail address}: \texttt{agustin.moreno2191@gmail.com}
 
}}

\date{}

\begin{document}

\maketitle

\begin{abstract}
    In this article, for Hamiltonian systems with two degrees of freedom, we study \emph{doubly symmetric} periodic orbits, i.e.\ those which are symmetric with respect to two (distinct) commuting antisymplectic involutions. These are ubiquitous in several problems of interest in mechanics. We show that, in dimension four, doubly symmetric periodic orbits cannot be negative hyperbolic. This has a number of consequences: (1) all covers of doubly symmetric orbits are \emph{good}, in the sense of Symplectic Field Theory \cite{eliashberg-givental-hofer}; (2) a non-degenerate doubly symmetric orbit is stable if and only if its CZ-index is odd; (3) a doubly symmetric orbit does \textbf{not} undergo period doubling bifurcation; and (4) there is always a stable orbit in any collection of doubly symmetric periodic orbits with negative \emph{SFT-Euler characteristic} (as coined in \cite{fkm}). The above results follow from: (5) a symmetric orbit is negative hyperbolic if and only its two \emph{$B$-signs} (introduced in \cite{frauenfelder-moreno}) differ.
\end{abstract}

\tableofcontents

\section{Introduction}

This article deals with the study of doubly symmetric periodic orbits in dimension four, i.e.\ for Hamiltonian systems with two degrees of freedom. These are ubiquitous in problems of interest in mechanics; we give several examples in Section~\ref{example}. Let us introduce the basic concepts.

\medskip

\begin{figure}
    \centering
    \includegraphics{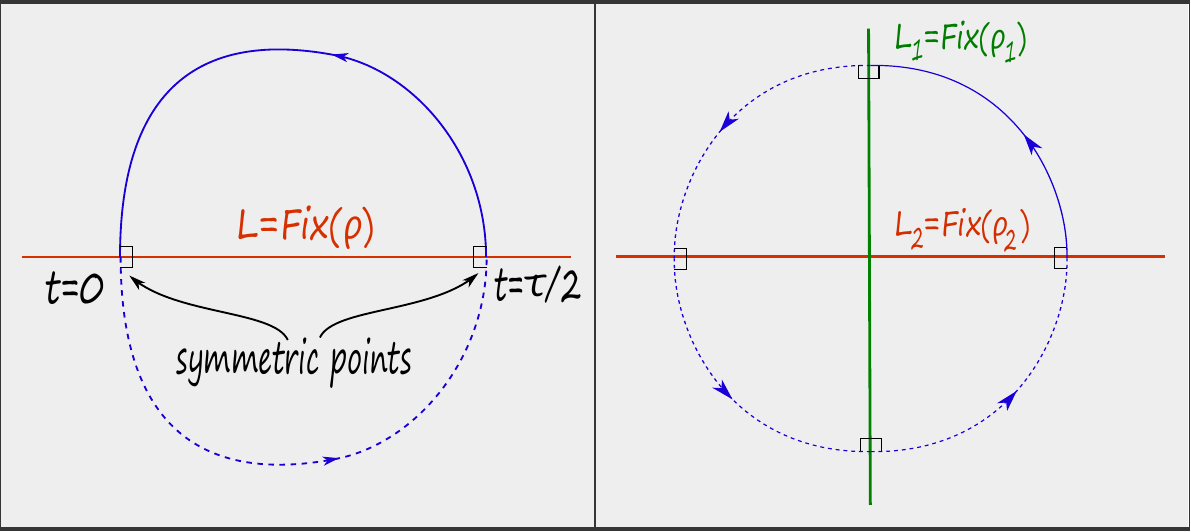}
    \caption{Left: A symmetric orbit. Right: A doubly symmetric orbit.}
    \label{fig:sym_orbit}
\end{figure}

\textbf{Symmetric orbits.} Consider a symplectic manifold $(M,\omega)$ endowed with an antisymplectic involution $\rho: M\rightarrow M$ (i.e.\ $\rho^2=id$, $\rho^*=\omega=-\omega)$, also referred to as a real structure. Its fixed point set $L=\mathrm{Fix}(\rho)$ is a Lagrangian submanifold of $M$. Given a Hamiltonian $H:M\rightarrow \mathbb R$, we say that $\rho$ is a \emph{symmetry} of the Hamiltonian system induced by $H$, if $H\circ \rho=H$. In this situation, a \emph{symmetric} periodic orbit is a periodic orbit $v: S^1=\mathbb R/\tau \mathbb R\rightarrow M$ satisfying $\rho(v(-t))=v(t)$ for all $t$. A symmetric periodic orbit can also be thought of as a chord starting and ending in $L$, where the endpoints coincide with $v(0), v(\tau/2)$ (the \emph{symmetric} points), see Figure \ref{fig:sym_orbit}.

Now suppose we have two \emph{distinct} antisymplectic involutions $\rho_1$ and $\rho_2$ which
commute with each other. In this case we have two Lagrangian submanifolds
$L_1=\mathrm{Fix}(\rho_1)$ and $L_2=\mathrm{Fix}(\rho_2)$. Given a chord from
$L_1$ to $L_2$ we can apply $\rho_2$ to it to get a chord from $L_1$ to itself.
Now apply $\rho_1$ to this chord. The resulting periodic orbit is then \emph{doubly} symmetric, as it is symmetric with respect to both $\rho_1,\rho_2$, see again Figure \ref{fig:sym_orbit}. We provide a more formal definition of
the notion of a doubly symmetric periodic orbits in Section~\ref{doublesym}.  

\medskip

\textbf{Reduced monodromy.} Suppose that $(M,\omega)$ is a four-dimensional symplectic manifold, $H\colon M \to \mathbb{R}$ is a smooth Hamiltonian, and $v$ is a nonconstant 
periodic orbit of the Hamiltonian vector field $X_H$ of $H$ of period $\tau$. By preservation of energy $H$ is constant along $v$, i.e., $v$ lies for
all times on a level set $\Sigma=H^{-1}(c)$ for some $c \in \mathbb{R}$. The differential of the flow $\phi^t_H$ induces a map on the two-dimensional
quotient vector space
$$M_v:=\overline{d\phi^\tau_H(v(0))} \colon T_{v(0)}\Sigma/\langle X_H(v(0))\rangle
\to T_{v(0)}\Sigma/\langle X_H(v(0))\rangle,$$
referred to as the \emph{reduced monodromy}. The reduced monodromy is a 
two-dimensional symplectic transformation, i.e., $\det M_v=1.$ Depending on the trace of its reduced monodromy, periodic orbits on a four-dimensional
symplectic manifold are now partitioned into three classes.
\begin{description}
 \item[Positive hyperbolic:] $\mathrm{tr}(M_v)>2,$
 in which case the reduced monodromy has two positive, real eigenvalues inverse to each other.
 \item[Negative hyperbolic:] $\mathrm{tr}(M_v)<2,$ in which case the reduced monodromy has two negative, real eigenvalues inverse to each other.
 \item[Elliptic:] $-2\leq \mathrm{tr}(M_v)\leq 2.$
 If the trace is precisely two, the reduced monodromy has one as an eigenvalue with
 algebraic multiplicity two. If the trace is precisely minus two, it has minus one
 as an eigenvalue with algebraic multiplicity two. Otherwise it has two nonreal
 eigenvalues on the unit circle conjugated to each other. 
\end{description}
In the language of Symplectic Field Theory, an even cover of a negative hyperbolic
orbit is called \emph{bad}; otherwise a periodic orbit is called \emph{good}. Here we prove the following:
\\ \\
\textbf{Theorem\,A: } \emph{For a Hamiltonian system with two degrees of freedom, a doubly symmetric periodic orbit cannot be negative hyperbolic.}
\\ \\
In particular, it follows from Theorem\,A that all covers of a doubly symmetric
periodic orbit are good periodic orbits. 

\medskip

\textbf{Stability.} While elliptic periodic orbits are stable, hyperbolic ones are unstable. 
On the other hand, elliptic and negative hyperbolic orbits have odd
Conley-Zehnder index, while positive hyperbolic ones have even Conley-Zehnder index.
For the second statement it is better to exclude the degenerate case where the trace of the reduced monodromy is two, since in this case there are different conventions on how to define 
the Conley-Zehnder index. We see from this that if we can exclude negative hyperbolic
orbits, the question of stability of a periodic orbit can be answered in terms
of the parity of its Conley-Zehnder index. In particular, we have the following
Corollary of Theorem\,A:
\\ \\
\textbf{Corollary\,B: } \emph{Suppose that $v$ is a nondegenerate doubly
symmetric periodic orbit of a Hamiltonian system with two degrees of freedom. Then it it stable if and only if its Conley-Zehnder
index is odd.}
\\ \\
\textbf{Overview of proof of Theorem\,A.} The proof of Theorem\,A uses a real
version of Krein theory for the reduced monodromy of a symmetric periodic orbit. Given a symmetric orbit $v$, the differential of the antisymplectic involution at $v(0) \in L=\mathrm{Fix}(\rho)$ induces an antisymplectic involution
$$R \colon T_{v(0)}\Sigma/\langle X_H(v(0))\rangle
\to T_{v(0)}\Sigma/\langle X_H(v(0))\rangle,$$
i.e.\ an orientation reversing involution on the two-dimensional vector space
$T_{v(0)}\Sigma/\langle X_H(v(0))$. 
The involution $R$ conjugates the reduced monodromy with its inverse, i.e.\
\begin{equation}\label{coninv}
RM_vR=M_v^{-1}.
\end{equation}
We choose a symplectic basis on $T_{v(0)}\Sigma/\langle X_H(v(0))$ such that
the involution $R$ gets identified with the matrix
$$R=\left(\begin{array}{cc}
1 & 0\\
0 &-1
\end{array}\right)$$
and the reduced monodromy is given by a matrix
$$M_v=\left(\begin{array}{cc}
a & b\\
c & d
\end{array}\right)$$
satisfying the determinant condition $ad-bc=1$. It follows from (\ref{coninv})
that $a=d$ so that 
$$M_v=\left(\begin{array}{cc}
a & b\\
c & a
\end{array}\right), \quad a^2-bc=1.$$
In particular, the question to which class the periodic orbit belongs is completely
answered by the entry $a$ of the reduced monodromy matrix. For fixed $a$, if an off-diagonal entry is not zero, then it completely determines the other one in view of the determinant
condition. On the other hand, the off-diagonal entries depends on the choice of the
symplectic basis used to identify the reduced monodromy with a matrix. Since
the symplectic basis vectors are required to be eigenvectors of the antisymplectic
involution $R$, such a symplectic basis is determined up to a scaling factor, so that
the identification of the reduced monodromy with a matrix is unique up to conjugation by
a matrix of the form
$$\left(\begin{array}{cc}
\mu & 0\\
0 & \frac{1}{\mu}
\end{array}\right), \quad  \mu \in \mathbb{R}\setminus \{0\}.$$
In particular, while the value of $b$ is not an invariant, its sign is an invariant.
Following \cite{frauenfelder-moreno} we refer to $\mathrm{sign}(b)$ as the
\emph{B-sign} of the reduced monodromy, see also \cite{zhou}. In the case elliptic case, by \cite[Appendix B]{frauenfelder-moreno}, the
B-sign gives the same information as the Krein type of the eigenvalues of the reduced monodromy (as introduced in \cite{kre1,kre2,K3,K4,Moser}). In the
hyperbolic case the eigenvalues have no Krein type. Therefore the B-sign in the hyperbolic
case is an additional invariant of the real structure $\rho$. 
\\ \\
A symmetric periodic orbit intersects the Lagrangian $L=\mathrm{Fix}(\rho)$ in its two
symmetric points. From the reduced monodromies of each symmetric point we obtain a B-sign,
so that a symmetric periodic orbit is actually endowed with two B-signs. The main
observation to prove Theorem\,A is the following:
\\ \\
\textbf{Theorem\,C: } \emph{A symmetric periodic orbit of a Hamiltonian system with two degrees of freedom is negative hyperbolic if and only
if its two B-signs are different.}
\\ \\
If the symmetric periodic orbit is elliptic it is actually clear that the two B-signs
have to agree. Indeed, as already mentioned, in the elliptic case the B-sign is just determined by the Krein sign of the eigenvalues. Since reduced monodromy matrices of a periodic orbit for different
starting points are all conjugated to each other, Theorem\,C follows in the elliptic case.
What remains to be examined is the hyperbolic case, namely
that in the positive hyperbolic case the two B-signs agree, while in the negative hyperbolic
case they disagree. To address this, in Section~\ref{realcouple} we introduce the notion
of \emph{real couples}, so that Theorem\,C becomes a consequence of Proposition~\ref{different} below. 
\\ \\
The strategy to prove Theorem\,A is now rather obvious. One shows that the additional
real structure for a doubly symmetric periodic orbit forces the two B-signs to agree,
so that, in view of Theorem\,C, a doubly periodic orbit cannot be negative hyperbolic. 
This is carried out in Section~\ref{redhol} where Theorem\,A is referred to as Corollary~\ref{notneghyp}.
\\ \\
\textbf{Period doubling bifurcation.} When considered in families, periodic orbits may undergo \emph{bifurcation}, by which a non-degenerate orbit becomes degenerate (i.e.\ $1$ becomes an eigenvalue of its monodromy), and new orbits may appear. \emph{Generic} bifurcations in dimension four are well understood, see e.g.\ \cite[p.\ 599]{abraham-marsden}. However, the presence of symmetry, and in particular the presence of doubly symmetric orbits, is non-generic, and hence one expects new phenomena. And indeed, what follows aligns well with this expectation.

As a particular case of bifurcations, the transition from an elliptic periodic orbit to a negative hyperbolic orbit leads to a \emph{period doubling bifurcation}, by which a new orbit appears, whose period is close to double the period of the original orbit. In the case where the negative hyperbolic orbit is symmetric, its two different B-signs can actually be useful to figure out where the new periodic orbit
of double period bifurcates, see \cite{frauenfelder-koh-moreno}. Namely, bifurcation happens near the symmetric point where the $B$-sign does \emph{not} jump. Moreover, a consequence of Theorem\,A is the following, which emphasizes the non-generic nature of symmetry:
\\ \\
\textbf{Corollary\,D: } \emph{In dimension four, doubly symmetric periodic orbits do \textbf{not} undergo period doubling bifurcation.}
\\ \\
Indeed, as in period doubling bifurcation the orbit itself does \emph{not} bifurcate (its double cover does), the orbit after such a bifurcation would have to be doubly symmetric if the orbit before bifurcation is, thus contradicting Theorem\,A. We remark that Corollary\,D fails in dimension six, i.e.\ for systems with three degrees of freedom. Indeed, see e.g.\ \cite[Section 6]{fkm} for a numerical example of a \emph{planar-to-spatial} period doubling bifurcation of doubly symmetric orbits.
\\ \\
\textbf{SFT-Euler characteristic.} In order to address the situation of more general bifurcations than period doubling bifurcation (in the presence of symmetry), we consider a \emph{Floer numerical invariant}. Namely, following \cite{frauenfelder-moreno}, the \emph{SFT-Euler characteristic} of a periodic orbit $v$ is by definition the Euler characteristic of its local Floer homology, given by
$$
\chi_{SFT}(v)=\#\{\mbox{good positive hyperbolic}\}-\#\{\mbox{elliptic, negative hyperbolic}\}.
$$
Here, one counts each type of orbit that appears after a generic perturbation of the orbit $v$, so that it bifurcates into a collection of non-degenerate orbits. We remark that bad orbits do not contribute to this number. Note also that this number is $\pm 1$ in the case where $v$ is itself non-degenerate, depending on its type. The remarkable fact, which follows from Floer theory, is that $\chi_{SFT}(v)$ is independent of the perturbation, and so in particular it remains invariant under bifurcations of $v$. It is therefore very useful in order to study non-generic bifurcations.

Moreover, given a collection of periodic orbits (which may not necessarily arise from a bifurcation, but e.g.\ as critical points of an action functional, with a priori fixed homotopy class) one can also consider the same number computed via the above formula. Its invariance under arbitrary homotopies will of course not be guaranteed, and will depend on the particular situation. An example of interest, for which a suitable homotopy invariance holds, are frozen planets. These are periodic orbits for the Helium problem which we discuss in more
detail in Section~\ref{example}. Due to the interaction between the two
electrons in Helium, frozen planets cannot be approached by perturbative methods but instead one can replace the instantaneous interaction
of the two electrons by a mean interaction. If one interpolates between
mean and instantaneous interaction one obtains a homotopy of a frozen planet problem for which one has compactness in the symmetric case
\cite{frauenfelder}. This allows one to define a version of the Euler characteristic for frozen planets which is invariant under this homotopy \cite{cieliebak-frauenfelder-volkov2}, and which agrees with the SFT-Euler charactersitic $\chi_{SFT}$ for the instantaneous interaction. The Euler characteristic for this problem is $-1$,
see the remark after Corollary\,B in \cite{cieliebak-frauenfelder-volkov2}. For each negative energy, this implies the existence of a symmetric frozen planet orbit for the instantaneous interaction, see Corollary C in \cite{cieliebak-frauenfelder-volkov}. This follows by homotopy invariance of the Euler characteristic, and the existence (proved analytically in \cite{f}) of a unique nondegenerate symmetric orbit for the mean interaction. 

With these motivations in mind, the following is again a consequence of Theorem\,A:

\medskip

\textbf{Corollary\,E: } \emph{In dimension four, suppose that a 
collection of doubly symmetric periodic orbits has negative SFT-Euler characteristic. Then a stable periodic orbit exists.}

\medskip

Indeed, Theorem\,A and the formula defining $\chi_{SFT}$ imply the existence of an elliptic orbit, and one needs to recall that elliptic orbits are precisely the stable orbits for a Hamiltonian system in dimension four.

\medskip

\textbf{Acknowledgements.} A.\ Moreno is supported by the National Science Foundation under Grant No.\ DMS-1926686, and by the Sonderforschungsbereich TRR 191 Symplectic Structures in Geometry, Algebra and Dynamics, funded by the DFG (Projektnummer 281071066 – TRR 191).

\section{Examples of doubly symmetric periodic orbits}\label{example}

\subsection{The direct and retrograde periodic orbit in Hill's lunar problem}

Hill's lunar Hamiltonian goes back to Hill's groundbreaking work on the orbit of
our Moon \cite{hill}, describing its motion around the Earth and the Sun. The Earth
lies in the center of the frame of reference, while the Sun, assumed to be infinitely much heavier than the
Earth, lies at infinity. The Hamiltonian reads
$$H \colon T^*(\mathbb{R}^2 \setminus \{0\}) \to \mathbb{R}, \quad
(q,p) \mapsto \frac{1}{2}\big((p_1+q_2)^2+(p_2-q_1)^2\big)
-\frac{1}{|q|}-\frac{3}{2}q_1^2.$$
It is invariant under the two commuting antisymplectic involutions
$$\rho_1, \rho_2 \colon T^*\mathbb{R}^2 \to T^*\mathbb{R}^2$$
given, for $(q,p) \in T^*\mathbb{R}^2$, by
$$\rho_1(q_1,q_2,p_1,p_2)=(q_1,-q_2,-p_1,p_2), \quad
\rho_2(q_1,q_2,p_1,p_2)=(-q_1,q_2,p_1,-p_2).$$
The fixed point sets of the two antisymplectic involutions are the conormal bundles of
the $x$-axis and the $y$-axis, respectively. If one studies a doubly symmetric periodic
orbit in configuration space $\mathbb{R}^2 \setminus \{0\}$, this means
that it starts perpendicularly at the $x$-axis, after a quarter period hits the
$y$-axis perpendicularly, then gets reflected at the $y$-axis for the next quarter
period, and finally gets reflected at the $x$-axis for the second half of the period. Such periodic orbits can be found by a shooting argument where one shoots perpendicularly
from the $x$-axis for a varying starting point at the $x$-axis, until one
hits the $y$-axis perpendicularly. Birkhoff used in \cite{birkhoff} this shooting 
argument to prove the existence of the retrograde periodic orbit for all energies
below the first critical value, see also \cite[Chapter 8.3.2]{frauenfelder-vankoert}.
Although the retrograde periodic orbit looks simpler than the direct one
\cite{henon}, astronomers are actually often more interested in the direct one, since
our Moon and actually most moons in our solar system are direct. However, there are
prominent counterexamples. Triton, the largest moon of the planet Neptun, is for example
retrograde.

\subsection{The Levi-Civita regularization}

Hill's lunar problem arises as a limit case of the restricted three-body problem, see
for instance \cite[Chapter 5.8.2]{frauenfelder-vankoert}. In the restricted three-body
problem the masses of the Sun and the Earth are comparable and their distance is finite.
Different from the Hill's lunar problem, the restricted three-body problem is only
invariant under the antisymplectic involution
$$\rho \colon T^* \mathbb{R}^2 \to T^* \mathbb{R}^2, \quad
(q_1,q_2,p_1,p_2) \mapsto (q_1,-q_2,-p_1,p_2)$$
obtained from reflection at the $x$-axis, but not anymore under the antisymplectic
involution corresponding to reflection at the $y$-axis. 
\\ \\
We identify $\mathbb{R}^2$ with the complex plane $\mathbb{C}$ and denote by $\mathbb{C}^*:=\mathbb{C} \setminus \{0\}$ the complex plane pointed at the origin. We consider the
squaring map
$$\ell \colon \mathbb{C}^* \to \mathbb{C}^*, \quad z \mapsto z^2.$$
Note that the squaring map
is a two-to-one covering. The contragradient (or symplectic lift) of the squaring map is the symplectic 
map
$$L \colon T^* \mathbb{C}^* \to T^*\mathbb{C}^*,
\quad (z,w) \mapsto \bigg(z^2, \frac{w}{2\bar{z}}\bigg),$$
where $\bar{z}$ is the complex conjugate of $z$.
This map was used by Levi-Civita to regularise two-body collisions 
\cite{levi-civita} and therefore it is known under the name of
\emph{Levi-Civita regularization}.
On $T^* \mathbb{C}$ we have the two commuting antisymplectic involutions
$$\sigma_1,\sigma_2 \colon T^* \mathbb{C} \to
T^*\mathbb{C}$$
which are given, for $(z,w) \in \mathbb{C} \times \mathbb{C}=T^*\mathbb{C}$, by
$$\sigma_1(z,w)=(\bar{z},-\bar{w}), \quad \sigma_2(z,w)=(-\bar{z},\bar{w}).$$
The Levi-Civita regularization lifts the restriction of the antisymplectic involution
$\rho$ to $T^* \mathbb{C}^*$ to the restriction of
$\sigma_1$ and $\sigma_2$ to $T^*\mathbb{C}^*$, so that we have
$$L \circ \sigma_1\big|_{T^*\mathbb{C}^*}=\rho\big|_{T^*\mathbb{C}^*} \circ L, \qquad
L \circ \sigma_2\big|_{T^*\mathbb{C}^*}=\rho\big|_{T^*\mathbb{C}^*} \circ L.$$
Now suppose that $v=(q,p)$ is a periodic orbit in $T^*\mathbb{C}^*$
which is symmetric with respect to $\rho$, and such that it has odd winding number around
the origin. Then $v$ lifts under the Levi-Civita regularisation to a periodic
orbit on $T^*\mathbb{C}^*$ which is doubly symmetric with respect
to $\sigma_1$ and $\sigma_2$. 

On the other hand, retrograde and direct orbits exist as
well in the restricted three-body problem. Different from Hill's lunar problem, they
are just symmetric, but not doubly symmetric. However, the lifts under the
Levi-Civita regularisation are doubly symmetric, as the retrograde and direct
periodic orbit have winding number one around the origin.

\subsection{Langmuir's periodic orbit}

\begin{figure}
    \centering
    \includegraphics{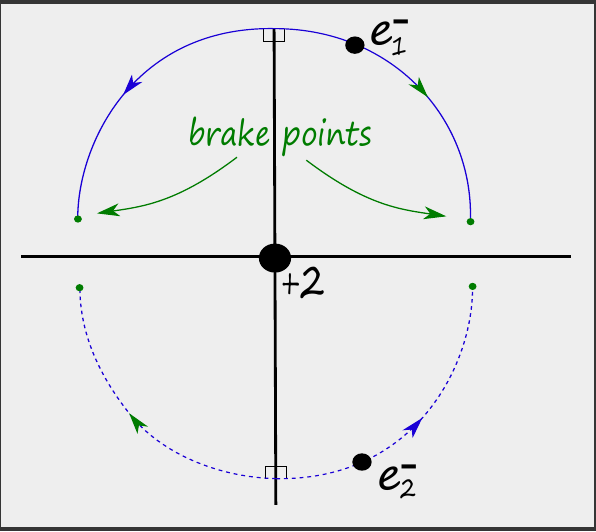}
    \caption{Langmuir's doubly symmetric orbit, and its symmetric version.}
    \label{fig:Langmuir}
\end{figure}

Langmuir's periodic orbit is a periodic orbit for the Helium problem. It was
first discovered by Langmuir \cite{langmuir} numerically as a candidate for the ground state of the Helium atom. For an analytic existence proof we refer to \cite{cieliebak-frauenfelder-schwingenheuer}, and for its role in the semiclassical treatment of Helium, to
\cite{tanner-richter-rost}.
\\ \\
In the Helium atom, there is a nucleus of positive charge plus two at the origin, i.e.\ there are two protons. It
attracts two electrons of charge minus one according to Coulomb's law, which looks
formally the same as Newton's law. Moreover, the two electrons repel each other, again according to Coulomb's law. We abbreviate by $$\Delta:=\big\{(q,q): q \in \mathbb{C}^*\big\} \subset \mathbb{C}^* \times
\mathbb{C}^*$$
the diagonal. The Hamiltonian for the \emph{planar} Helium problem is then a smooth function
$$H \colon T^*\big(\mathbb{C}^* \times \mathbb{C}^* \setminus \Delta\big) \to 
\mathbb{R}$$
given by
$$H(q_1,q_2,p_1,p_2)=\frac{1}{2}|p_1|^2+\frac{1}{2}|p_2|^2-\frac{2}{|q_1|}
-\frac{2}{|q_2|}+\frac{1}{|q_1-q_2|}.$$
The Hamiltonian is invariant under the symplectic involution
$$\sigma \colon T^*\big(\mathbb{C}^* \times \mathbb{C}^* \setminus \Delta\big)
\to T^*\big(\mathbb{C}^* \times \mathbb{C}^* \setminus \Delta\big)$$
given by
$$\sigma(q_1,q_2,p_1,p_2)=(\bar{q}_2,\bar{q}_1,\bar{p}_2,\bar{p}_1),$$
consisting of the combination of particle interchange and reflection at the
$x$-axis. The \emph{Langmuir Hamiltonian} is the restriction of $H$ to
the fixed point set of $\sigma$
$$H_\sigma:=H\big|_{\mathrm{Fix}(\sigma)} \colon \mathrm{Fix}(\sigma) \to \mathbb{R}.$$
The fixed points set consists of points $(q_1,q_2,p_1,p_2)
\in T^*(\mathbb{C}^* \times \mathbb{C}^* \setminus \Delta)$ which satisfy
$$q_1=\bar{q}_2=:q, \quad p_1=\bar{p}_2=:p.$$
It therefore suffices to consider the Langmuir Hamiltonian on the cotangent
bundle of the upper halfplane
$$\mathbb{H}=\big\{q=q_1+iq_2 \in \mathbb{C}: q_2>0\big\}$$
where it is given by
$$H_\sigma(q,p)=|p|^2-\frac{4}{|q|}+\frac{1}{2q_2}.$$
On the cotangent bundle of the uper halfplane we have the two antisymplectic involutions
$$\rho_1,\rho_2 \colon T^* \mathbb{H} \to T^*\mathbb{H}$$
given by
$$\rho_1(q,p)=(-\bar{q},\bar{p}), \quad \rho_2(q,p)=(q,-p),$$
under both of which $H_\sigma$ is invariant. The fixed point set of $\rho_1$
is the conormal bundle of the positive imaginary axis, while the fixed point
set of $\rho_2$ consists of \emph{brake points}, i.e.\ at which the velocity is zero. The Langmuir orbit for the first electron $e_1^-$ starts perpendicularly at the imaginary axis and brakes at a quarter of the period,
and is therefore a doubly symmetric periodic orbit with respect to $\rho_1$ and 
$\rho_2$. The second electron $e_2^-$ similarly has an associated Langmuir orbit, obtained by conjugation of that of $e_1^-$, see Figure \ref{fig:Langmuir}.

\subsection{Symmetric frozen planets}

\begin{figure}
    \centering
    \includegraphics{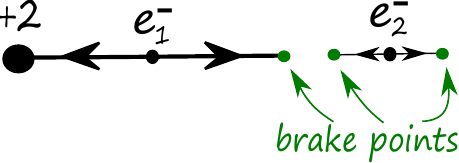}
    \caption{A frozen planet configuration.}
    \label{fig:frozen_planet}
\end{figure}

Other examples of periodic orbits for the Helium problem are frozen planet orbits.
In this examples both electrons lie on a line on the same side of the nucleus.
The inner electron makes consecutive collisions with the nucleus. The outer electron,
the actual ``frozen planet", which is attracted by the nucleus but repelled by the
inner electron, stays almost stationary but librates slightly. Frozen planet orbits were discovered by
physicists \cite{tanner-richter-rost, wintgen-richter-tanner} in the context of
semiclassics. They recently attracted the interest of mathematicians
\cite{cieliebak-frauenfelder-volkov, zhao}. A frozen planet orbit is
called \emph{symmetric} if the two electrons brake at the same time, and at the time
the inner electron collides with the nucleus the outer electron brakes again, see Figure \ref{fig:frozen_planet}. If one applies the Levi--Civita regularization to a symmetric frozen planet one obtains a doubly symmetric periodic orbit. 

\section{Real couples}\label{realcouple}

A \emph{real symplectic vector space} is a triple $(V,\omega,R)$ consisting
of a symplectic vector space $(V,\omega)$ and a linear antisymplectic involution
$R \colon V \to V$, i.e.\ $R^2=Id, R^*\omega=-\omega$. 
\begin{fed}
Assume that $(V_1,\omega_1,R_1)$ and $(V_2,\omega_2,R_2)$ are real
symplectic vector spaces. A \textbf{real couple} $(\Psi,\Phi)$ 
is a tuple of linear symplectic maps
$$\Psi \colon (V_1,\omega_1) \to (V_2,\omega_2), \quad \Phi \colon (V_2,\omega_2) \to (V_1,\omega_1)$$
which are related by
\begin{equation}\label{couple}
R_2 \Psi R_1=\Phi^{-1}.
\end{equation}
\end{fed}
Note that if $(\Psi,\Phi)$ is a real couple, then $(\Phi,\Psi)$
is one as well, since it follows from (\ref{couple}) that
$$R_1 \Phi R_2=R_1 R_1^{-1} \Psi^{-1} R_2^{-1}R_2=\Psi^{-1}.$$
If $(\Psi,\Phi)$ is a real couple then its composition
$$\Phi \Psi \colon (V_1,\omega_1) \to (V_1,\omega_1)$$
is a linear symplectic map from the fixed symplectic vector space
$(V_1,\omega_1)$ into itself which has the special property that it
is conjugated to its inverse via the antisymplectic involution $R_1$. Indeed,
\begin{eqnarray}\label{conj}
R_1 \Phi \Psi R_1&=&R_1 \Phi R_2 R_2 \Psi R_1=\Psi^{-1} \Phi^{-1}=(\Phi \Psi)^{-1}.
\end{eqnarray}
We now consider more closely the two-dimensional case. Note that every two-dimensional
real symplectic vector space is conjugated to $\mathbb{R}^2$, endowed with its
standard symplectic structure and antisymplectic involution
$$R=\left(\begin{array}{cc}
1 & 0\\
0 & -1
\end{array}\right).$$
After such conjugation, a real couple then consists of a pair of matrices
$$(A,B) \in SL(2;\mathbb{R}) \times SL(2;\mathbb{R})$$
such that
\begin{equation}\label{conjmat}
RAR=B^{-1}.
\end{equation}
Writing 
$$A=\left(\begin{array}{cc}
a & b\\
c & d
\end{array}\right), \quad ad-bc=1$$
we have
\begin{eqnarray*}
RAR&=&\left(\begin{array}{cc}
1 & 0\\
0 & -1
\end{array}\right)\left(\begin{array}{cc}
a & b\\
c & d
\end{array}\right)
\left(\begin{array}{cc}
1 & 0\\
0 & -1
\end{array}\right)\\
&=&\left(\begin{array}{cc}
1 & 0\\
0 & -1
\end{array}\right)\left(\begin{array}{cc}
a & -b\\
c & -d
\end{array}\right)\\
&=&\left(\begin{array}{cc}
a & -b\\
-c & d
\end{array}\right)
\end{eqnarray*}
and therefore
$$B=(RAR)^{-1}=\left(\begin{array}{cc}
d & b\\
c & a
\end{array}\right).$$
Hence their products are given by the following matrices
\begin{equation}\label{prod1}
AB=\left(\begin{array}{cc}
a & b\\
c & d
\end{array}\right)\left(\begin{array}{cc}
d & b\\
c & a
\end{array}\right)=\left(\begin{array}{cc}
ad+bc & 2ab\\
2cd & ad+bc
\end{array}\right)
\end{equation}
and
\begin{equation}\label{prod2}
BA=\left(\begin{array}{cc}
d & b\\
c & a
\end{array}\right)\left(\begin{array}{cc}
a & b\\
c & d
\end{array}\right)=\left(\begin{array}{cc}
ad+bc & 2bd\\
2ac & ad+bc
\end{array}\right).
\end{equation}
Since
$$BA=B(AB)B^{-1}$$
the two product are conjugated to each other in $SL(2;\mathbb{R})$. Moreover, 
they both belong to the subspace 
$$SL^R(2;\mathbb{R}):=\Bigg\{M=\left(\begin{array}{cc}
\alpha & \beta \\
\gamma & \alpha
\end{array}\right): \alpha^2-\beta \gamma=1\Bigg\}$$
of $SL(2;\mathbb{R})$. If
$M \in SL^R(2;\mathbb{R})$ satisfies $\mathrm{tr}(M) \neq \pm 2$ we 
define its \emph{real Krein sign} as
$$\kappa(M):=\mathrm{sign}(\beta).$$
Note that the trace condition implies that $\alpha \neq \pm 1$ so that, in view
of the determinant condition $\alpha^2-\beta \gamma$, we have that $\beta \neq 0$,
and so its sign is well defined. The following proposition is now straightforward to
prove.
\begin{prop}\label{different}
The real Krein signs of $AB$ and $BA$ differ, if and only if
\begin{equation}\label{tra}
\mathrm{tr}(AB)=\mathrm{tr}(BA)<-2,
\end{equation}
i.e., if and only if $AB$ and therefore as well $BA$ are negative hyperbolic. 
\end{prop}
\textbf{Proof: } By (\ref{prod1}) and (\ref{prod2}) the trace condition (\ref{tra})
is equivalent to the inequality
$$ad+bc<-1.$$
In view of the determinant condition $ad-bc=1$ this in turn is equivalent to the 
inequality 
$$ad<0,$$
i.e., the requirement that the signs of $a$ and $d$ are different. Having once more
a look at (\ref{prod1}) and (\ref{prod2}), we see that this happens if and only if
the real Krein signs of $AB$ and $BA$ disagree. This proves the proposition. \hfill $\square$
\\ \\
In the following we assume that $(\Psi,\Phi)$ is a real couple between real
symplectic vector spaces $(V_1,\omega_1,R_1)$ and $(V_2,\omega_2,R_2)$.
\begin{fed}
The real couple $(\Psi,\Phi)$ is called \textbf{symmetric} if there exists
a linear map 
$$S \colon V_1 \to V_2$$
which is antisymplectic, i.e.,
$$S^* \omega_2=-\omega_1$$
and satisfies
\begin{equation}\label{symcoup}
\Psi=S\Psi^{-1}S, \quad \Phi^{-1}=S\Phi S, \quad R_2 S R_1=S.
\end{equation}
\end{fed}
For a symmetric real couple
$$T:=SR_1=R_2 S \colon (V_1,\omega_1) \to (V_2,\omega_2)$$
is a linear symplectic map which in view of
$$TR_1=S=R_2 T$$
interchanges the two real structures, so that $T$ leads to an identification of
the two real symplectic vector spaces $(V_1,\omega_1,R_1)$ and $(V_2,\omega_2,R_2)$. 
In the two-dimensional case if we identify this further with $\mathbb{R}^2$
endowed with its standard symplectic form and standard real structure $R$, then
not only $R_1$ and $R_2$ are identified with $R$, but so is $S$. The real tuple
becomes identified with a pair $(A,B)$ of $SL(2,\mathbb{R})$-matrices which not
only satisfy (\ref{conjmat}) but due to (\ref{symcoup}) also satisfy
$$RAR=A^{-1}, \qquad RBR=B^{-1},$$
i.e., both matrices are conjugated to their inverse via $R$ and therefore lie
in the subspace $SL^R(2;\mathbb{R})$ of $SL(2;\mathbb{R})$. This implies that
$$A=B=\left(\begin{array}{cc}
a & b\\
c & a
\end{array}\right),\quad a^2-bc=1$$
and therefore
$$AB=BA.$$
In particular, $AB$ and $BA$ have the same real Krein sign. Therefore we obtain the
following corollary from Proposition \ref{different}.
\begin{cor}\label{notbad}
Suppose that $(\Psi,\Phi)$ is a two-dimensional symmetric real couple. Then neither
$\Phi \Psi$ nor $\Psi \Phi$ are negative hyperbolic. 
\end{cor}

\section{Doubly symmetric periodic orbits}\label{doublesym}

Suppose that $(M,\omega)$ is a symplectic manifold and $H \colon M \to \mathbb{R}$
is a smooth Hamiltonian. The Hamiltonian vector field $X_H$ of $H$ is implicitly defined
by the condition
$$dH=\omega(\cdot, X_H).$$
We abbreviate by $S^1=\mathbb{R}/\mathbb{Z}$ the circle. A \emph{simple periodic orbit} is
a bijective map $v \colon S^1 \to \mathbb{R}$ for which there exists $\tau>0$ such
that $v$ solves the ODE
$$\partial_t v(t)=\tau X_H(v(t)), \quad t \in S^1.$$
Since for a simple periodic orbit the map is bijective the Hamiltonian vector field
$X_H$ is nonvanishing along $v$ and therefore $\tau$ is uniquely determined by $v$.
We refer to $\tau$ as the \emph{period} of the simple periodic orbit $v$. We abbreviate
by
$$\mathcal{P}_H \subset C^\infty(S^1,M)$$
the set of simple periodic orbits of the Hamiltonian vector field $X_H$. 
\\ \\
A \emph{real symplectic manifold} is a triple $(M,\omega,\rho)$ where
$(M,\omega)$ is a symplectic manifold and $\rho \in \mathrm{Diff}(M)$ is an antisymplectic
involution on $M$, i.e.,
$$\rho^2=\mathrm{id}, \quad \rho^* \omega=-\omega.$$
If $H \colon M \to \mathbb{R}$ is a smooth function on a real symplectic manifold
which is invariant under the antisymplectic involution, i.e., 
$$H \circ \rho=H,$$
then its Hamiltonian vector field is anti-invariant, i.e.,
$$\rho^*X_H=-X_H.$$
We then obtain an involution 
$$I \colon \mathcal{P}_H \to \mathcal{P}_H, \quad v \mapsto \rho \circ v^-$$
where $v^-$ is the orbit traversed backwards, i.e.,
$$v^-(t)=v(-t), \quad t \in S^1.$$
A \emph{simple symmetric periodic orbit} is a fixed point of $I$, i.e,
$v \in \mathcal{P}_H$ satisfying
$$I(v)=v.$$
We abbreviate by
$$\mathcal{P}^I_H:=\mathrm{Fix}(I) \subset \mathcal{P}_H$$
the set of simple symmetric periodic orbits. 
We remark that the fixed point set of an antisymplectic involution
$$L:=\mathrm{Fix}(\rho)$$
is a Lagrangian submanifold of $M$. Note that if $v \in \mathcal{P}^I_H$ then
$$v\big(0\big), v\big(\tfrac{1}{2}\big) \in L$$
so that $v_{[0,\frac{1}{2}]}$ can be interpreted as a chord from $L$ to $L$. 
\\ \\
A \emph{doubly real symplectic manifold} is a quadruple $(M,\omega,\rho_1,\rho_2)$
where $(M,\omega)$ is a symplectic manifold and $\rho_1, \rho_2 \in \mathrm{Diff}(M)$
are two distinct antisymplectic involutions which commute with each other. Note since
$\rho_1$ and $\rho_2$ commute their composition
$$\sigma:=\rho_1 \circ \rho_2=\rho_2 \circ \rho_1$$
is a symplectic involution on $(M,\omega)$. 
Suppose that $(M,\omega,\rho_1,\rho_2)$ is a doubly real symplectic manifold
and $H \colon M \to \mathbb{R}$ is a smooth map which is invariant under both
involutions $\rho_1$ and $\rho_2$. We then have on the set of
simple periodic orbits $\mathcal{P}_H$ two involutions 
$$I_1 \colon \mathcal{P}_H \to \mathcal{P}_H, \quad v \mapsto \rho_1 \circ v^-, 
\qquad I_2 \colon \mathcal{P}_H \to \mathcal{P}_H, \quad v \mapsto \rho_2 \circ v^-.$$
Moreover, we have two Lagrangian submanifolds of $M$
$$L_1=\mathrm{Fix}(\rho_1), \qquad L_2=\mathrm{Fix}(\rho_2).$$
\begin{fed}
Suppose that $(M,\omega,\rho_1,\rho_2)$ is a doubly real symplectic manifold and
$H \colon M \to \mathbb{R}$ is a smooth function invariant under both involutions
$\rho_1$ and $\rho_2$.
A simple symmetric periodic orbit $v \in \mathcal{P}^{I_1}_H$ of $\rho_1$ is called
\textbf{doubly symmetric} if
\begin{equation}\label{dsym}
\rho_2 \circ v\big(0\big)=v\big(\tfrac{1}{2}\big).
\end{equation}
\end{fed}
Observe that since for a symmetric periodic orbit $v(1/2)$ lies in the fixed point set of
$\rho_1$ condition (\ref{dsym}) is equivalent to
$$\sigma \circ v\big(0\big)=v\big(\tfrac{1}{2}\big).$$
Doubly symmetric periodic orbits with respect to $\rho_1$ are in natural one-to-one
correspondence with double symmetric periodic orbits with respect to $\rho_2$. For
$r \in S^1$ and $v \in \mathcal{P}_H$ we denote by
$$r_*v \in \mathcal{P}_H$$
the reparametrized simple periodic orbit
$$r_*v (t)=v(r+t), \quad t \in S^1.$$
We have the following lemma.
\begin{lemma}\label{dore}
An orbit $v \in \mathcal{P}_H^{I_1}$ is doubly symmetric with respect to $\rho_2$
if and only if $\big(\tfrac{1}{4}\big)_*v \in \mathcal{P}_H^{I_2}$ is doubly symmetric
with respect to $\rho_1$. 
\end{lemma}
\textbf{Proof: }Suppose that $v \in \mathcal{P}_H^{I_1}$ is doubly symmetric with
respect to $\rho_2$. After reparametrization a simple periodic orbit is still a
simple periodic orbit so that we have
$$\big(\tfrac{1}{4}\big)_* v \in \mathcal{P}_H.$$
Since $H$ is invariant under $\rho_2$ we have that
$$I_2 \Big(\big(\tfrac{1}{4}\big)_* v\Big) \in \mathcal{P}_H.$$
Using (\ref{dsym}) we compute
\begin{eqnarray*}
I_2 \Big(\big(\tfrac{1}{4}\big)_* v\Big)\big(\tfrac{1}{4}\big)&=&
\rho_2 \circ \Big(\big(\tfrac{1}{4}\big)_* v\Big)^-\big(\tfrac{1}{4}\big)\\
&=&\rho_2\Big(\big(\tfrac{1}{4}\big)_*v\Big)\big(-\tfrac{1}{4}\big)\\
&=&\rho_2 \circ v\big(\tfrac{1}{4}-\tfrac{1}{4}\big)\\
&=&\rho_2 \circ v(0)\\
&=&v\big(\tfrac{1}{2}\big)\\
&=&\Big(\big(\tfrac{1}{4}\big)_*v\Big)\big(\tfrac{1}{4}\big).
\end{eqnarray*}
That means that $\big(\tfrac{1}{4}\big)_*v$ and $I_2\big(\big(\tfrac{1}{4}\big)_*v\big)$
are solutions of the same first order ODE which at time $\tfrac{1}{4}$ go through the
same point. Therefore from the uniqueness of the initial value problem of first order
ODE's we deduce that
$$I_2\Big(\big(\tfrac{1}{4}\big)_*v\Big)=\big(\tfrac{1}{4}\big)_*v$$
and hence
$$\big(\tfrac{1}{4}\big)_*v \in \mathcal{P}^{I_2}_H.$$
It remains to check its double symmetry with respect to $\rho_1$. For that we compute
\begin{eqnarray*}
\rho_1 \circ \Big(\big(\tfrac{1}{4}\big)_* v\Big)(0)&=&
\rho_1 \circ v\big(\tfrac{1}{4}\big)\\
&=&v\big(-\frac{1}{4}\big)\\
&=&v\big(\tfrac{3}{4}\big)\\
&=&\Big(\big(\tfrac{1}{4}\big)_*v\Big)\big(\tfrac{1}{2}\big).
\end{eqnarray*}
Here we have used in the second equation that $v$ is symmetric with respect to
$\rho_1$ and in the third equation that it is one-periodic. This shows that
$\big(\frac{1}{4}\big)_* v$ is doubly symmetric with respect to $\rho_1$.
\\ \\
It remains to check that if $\big(\frac{1}{4}\big)_* v \in \mathcal{P}^{I_2}_H$ is doubly symmetric with
respect to $\rho_1$ it follows that $v \in \mathcal{P}^{I_1}_H$ is doubly symmetriy with respect to $\rho_2$. Interchanging in the previous discussion the roles of $\rho_1$ and
$\rho_2$ we obtain that
$$\big(\tfrac{1}{4}\big)_*\big(\tfrac{1}{4}\big)_*v=\big(\tfrac{1}{2}\big)_*v
\in \mathcal{P}^{I_1}_H$$
is doubly symmetric with respect to $\rho_2$. The fact that $\big(\tfrac{1}{2}\big)_*v$
is invariant under $I_1$ implies that
\begin{eqnarray*}
I_1 v(t)&=&\rho_1 \circ v^-(t)\\
&=& \rho_1 \circ v(-t)\\
&=& \rho_1 \circ \Big(\big(\tfrac{1}{2}\big)_* v\Big)\big(-t-\tfrac{1}{2}\big)\\
&=& \rho_1 \circ \Big(\big(\tfrac{1}{2}\big)_* v\Big)^-\big(t+\tfrac{1}{2}\big)\\
&=& I_1\Big(\big(\tfrac{1}{2}\big)_* v\Big)\big(t+\tfrac{1}{2}\big)\\
&=&\Big(\big(\tfrac{1}{2}\big)_*v\Big)\big(t+\tfrac{1}{2}\big)\\
&=&v\big(t+1)\\
&=&v(t),
\end{eqnarray*}
so that $v \in \mathcal{P}^{I_1}_H$ is as well invariant under $I_1$. Since
$\big(\tfrac{1}{2}\big)_*v$ is doubly symmetric with respect to $\rho_2$ we obtain further
that
\begin{eqnarray*}
\rho_2 \circ v(0)&=&\rho_2 \circ \Big(\big(\tfrac{1}{2}\big)_* v\Big)\big(-\tfrac{1}{2}\big)\\
&=&\rho_2 \circ \Big(\big(\tfrac{1}{2}\big)_* v\Big)\big(\tfrac{1}{2}\big)\\
&=&\rho_2^2 \circ \Big(\big(\tfrac{1}{2}\big)_* v\Big)\big(0\big)\\
&=&\Big(\big(\tfrac{1}{2}\big)_* v\Big)\big(0\big)\\
&=&v\big(\tfrac{1}{2}\big),
\end{eqnarray*}
so that $v$ is doubly symmetric with respect to $\rho_2$ as well. This finishes
the proof of the lemma. \hfill $\square$

\section{The reduced monodromy}\label{redhol}

Suppose that $(M,\omega)$ is a symplectic manifold and 
$H \colon M \to \mathbb{R}$ is a smooth function. We denote by $\phi^t_H$
the flow of the Hamiltonian vector field of $H$, characterized by 
$$\phi^0_H(x)=x,\quad \frac{d}{dt}\phi^t_H(x)=X_H(\phi^t_H(x)), \qquad x \in M.$$
If $v$ is a simple periodic orbit of $X_H$ of period $\tau$ we have
$$\phi^\tau_H(v(0))=v(0),$$
i.e., $v(0)$ is a fixed point of $\phi^\tau_H$. The differential of the flow
$$d\phi^\tau_H(v(0)) \colon T_{v(0)} M \to T_{v(0)} M$$
is a linear symplectic map of the symplectic vector space
$(T_{v(0)}M,\omega_{v(0)})$ into itself. This map is referred to as
the \emph{unreduced monodromy}. Since $H$ is autonomous, i.e., does not
depend on time, we have 
$$d\phi^\tau_H(v(0))X_H(v(0))=X_H(v(0)).$$
Moreover, by preservation of energy the Hamiltonian $H$ is preserved along the flow
of its Hamiltonian vector field. In particular, if $c$ is the energy of $v$, i.e., the value
$H$ attains along $v$, the differential of the flow maps the 
tangent space $T_{v(0)}\Sigma$ of the energy hypersurface 
$$\Sigma=H^{-1}(c)$$
back to itself. Therefore
the unreduced monodromy induces a linear map
$$M_v:=\overline{d \phi^\tau_H(v(0))} \colon T_{v(0)}\Sigma/\langle X_H(v(0))\rangle
\to T_{v(0)}\Sigma/\langle X_H(v(0))\rangle$$
which is still symplectic for the symplectic structure on
$T_{v(0)}\Sigma/\langle X_H(v(0))\rangle$ induced from $\omega_{v(0)}$. This map
is referred to as the \emph{reduced monodromy}. Instead of restricting our attention to $0$ we could consider the reduced monodromy 
$$M_v^t:=\overline{d \phi^\tau_H(v(t))} \colon T_{v(t)}\Sigma/\langle X_H(v(t))\rangle
\to T_{v(t)}\Sigma/\langle X_H(v(t))\rangle$$
for any $t \in S^1$. Note that for different times $t$ the reduced monodromies are
symplectically conjugated to each other by the flow. 
\\ \\
Suppose now in addition that $\rho$ is a real structure on $(M,\omega)$ under which
$H$ is invariant and
$v \in \mathcal{P}_H^I$ is a symmetric periodic orbit. Since both points
$v(0)$ and $v\big(\tfrac{1}{2}\big)$ lie in the fixed point set 
of $\rho$ the differential of $\rho$ gives rise to linear antisymplectic involutions
$$d\rho\big(v\big(0\big)\big)  \colon T_{v(0)} M \to T_{v(0)} M, \quad
d\rho\big(v\big(\tfrac{1}{2}\big)\big) \colon T_{v(1/2)} M \to T_{v(1/2)} M
$$
which induce real structures on the quotient spaces
$T_{v(0)}\Sigma/\langle X_H(v(0))\rangle$ respectively 
$T_{v(1/2)}\Sigma /\langle X_H(v(1/2))\rangle$.
Since
the Hamiltonian vector field is anti-invariant, the antisymplectic involution $\rho$
conjugates the forward flow to the backward flow
$$\rho \phi^t_H \rho=\phi^{-t}_H.$$
In particular, differentiating this identity we have
$$d\rho\big(v\big(\tfrac{1}{2}\big)\big) \circ
d\phi^{\tau/2}_H\big(v\big(0\big)\big) \circ
d\rho\big(v\big(0\big)\big)=\Big(d \phi_H^{\tau/2}\big(v\big(\tfrac{1}{2}\big)
\big)\Big)^{-1}.$$
Therefore the induced maps
$$\Psi:=\overline{d\phi^{\tau/2}_H\big(v\big(0\big)\big)}
\colon T_{v(0)}\Sigma/\langle X_H(v(0))\rangle
\to T_{v(1/2)}\Sigma/\langle X_H(v(1/2))\rangle$$
and 
$$\Phi:=\overline{d\phi^{\tau/2}_H\big(v\big(\tfrac{1}{2}\big)\big)}
\colon T_{v(1/2)}\Sigma/\langle X_H(v(1/2))\rangle
\to T_{v(0)}\Sigma/\langle X_H(v(0))\rangle$$
give rise to a real couple $(\Psi,\Phi)$. Note that the compositions
coincide with the reduced monodromies at times $0$ and $\tfrac{1}{2}$ 
$$\Phi \Psi=\overline{d\phi^{\tau}_H\big(v\big(0\big)\big)}, \quad
\Psi \Phi=\overline{d\phi^{\tau}_H\big(v\big(\tfrac{1}{2}\big)\big)}.$$
Now we even assume that the symplectic manifold $(M,\omega)$ is doubly real
with real structures $\rho_1$ and $\rho_2$ under both of which 
$H$ is invariant and $v \in \mathcal{P}^{I_1}_H$ is doubly symmetric with
respect to $\rho_2$. The differential of $\rho_2$ gives rise to a linear antisymplectic
map
$$d\rho_2(v(0)) \colon T_{v(0)}M \to T_{v(1/2)}M$$
which induces an antisymplectic map on the quotient spaces
$$S \colon T_{v(0)}\Sigma/\langle X_H(v(0))\rangle \to T_{v(1/2)}\Sigma
/\langle X_H(v(1/2))\rangle.$$
Since $\rho_1$ commutes with $\rho_2$ this map interchanges the real structures.
By Lemma~\ref{dore} we have that $\big(\tfrac{1}{4}\big)_*v \in 
\mathcal{P}^{I_2}_H$ and therefore $S$ makes the real couple $(\Psi,\Phi)$
symmetric. Therefore we obtain, as a consequence of Corollary~\ref{notbad}, the following
corollary, which is Theorem\,A from the Introduction:
\begin{cor}\label{notneghyp}
A doubly symmetric periodic orbit on a four-dimensional symplectic manifold cannot
be negative hyperbolic. 
\end{cor}

\Addresses


\begin{thebibliography}{99}
\bibitem{abraham-marsden}R.\,Abraham, J.\,Marsden,
 \emph{Foundations of Mechanics}, 2nd ed. Addison-Wesley, New York
 (1978).
\bibitem{birkhoff} G.\,Birkhoff, \emph{The restricted problem of
 three bodies}, Rend.\,Circ.\,Matem.\,Palermo \textbf{39} (1915),
 265--334.
\bibitem{cieliebak-frauenfelder-schwingenheuer} K.\,Cieliebak, U.\,Frauenfelder,
M.\,Schwingenheuer, \emph{On Langmuir's periodic orbit}, Arch.\,Math. (Basel)
\textbf{118} (2022), no.\,4, 413--425.
\bibitem{cieliebak-frauenfelder-volkov} K.\,Cieliebak, U.\,Frauenfelder, E.\,Volkov,
\emph{A variational approach to frozen planet orbits in helium}, to appear in
Ann.\,Inst.\,H.\,Poincar\'e.
\bibitem{cieliebak-frauenfelder-volkov2}
K.\,Cieliebak, U.\,Frauenfelder, E.\,Volkov,
\emph{Nondegeneracy and integral count of
frozen planets in Helium}, arXiv: 2209.12634
\bibitem{eliashberg-givental-hofer} Y.\,Eliashberg, A.\,Givental, H.\,Hofer \emph{Introduction to Symplectic Field Theory}, Geom.\,Funct.\,Anal. 2000, Special Volume,
Part II, 560--673.

\bibitem{f} U.\,Frauenfelder, \emph{Helium and Hamiltonian delay equations}, Israel Journal of
Mathematics 246, 239–260 (2021).

\bibitem{frauenfelder} U.\,Frauenfelder, \emph{A compactness theorem
for frozen planets}, arXiv: 2010:15532, to appear in J.\,Topology and
Analysis.
\bibitem{frauenfelder-koh-moreno} U.\,Frauenfelder, D.\,Koh, U.\,Frauenfelder,
\emph{On Floer-type numerical invariants, GIT quotients, and orbit bifurcations
 of real-life planetary systems}, arXiv:2206.00627
\bibitem{frauenfelder-moreno} U.\,Frauenfelder, A.\,Moreno, \emph{On GIT quotients
of the symplectic group, stability and bifurcations of symmetric orbits}, arXiv:2109.09147

\bibitem{fkm} U.\ Frauenfelder, D.\ Koh, A.\ Moreno, \emph{Symplectic methods in the numerical search of orbits in real-life planetary systems
}, Preprint	arXiv:2206.00627. 

\bibitem{frauenfelder-vankoert} U.\,Frauenfelder, O.\,van Koert,\emph{The restricted
three-body problem and holomorphic curves}, Pathways in Mathematics, Birkh\"auser/Springer,
Cham (2018).
\bibitem{henon}M.\,H\'enon, \emph{Numerical exploration of the restricted problem. V.  Hill's case: Periodic orbits and their stability}, Astron.\,Astrophys. \textbf{1}, 223--238.
\bibitem{hill}G.\,Hill, \emph{Researches in the lunar theory}, Amer.\,J.\,Math.\textbf{1}
(1878), 5--26, 129--147, 245--260.
\bibitem{kre1} Krein, M.: Generalization of certain investigations of A.M. Liapunov on linear
differential equations with periodic coefficients. Doklady Akad. Nauk USSR 73 (1950) 445-448.

\bibitem{kre2} M.\ Krein, \emph{On the application of an algebraic proposition in the theory of monodromy matrices}, Uspekhi Math. Nauk 6 (1951) 171-177.

\bibitem{K3} M.\ Krein, \emph{On the theory of entire matrix-functions of exponential type,} Ukrainian Math. Journal 3 (1951) 164-173.

\bibitem{K4} M.\ Krein, \emph{On some maximum and minimum problems for characteristic numbers and Liapunov stability zones.,} Prikl. Math. Mekh. 15 (1951) 323-348.

\bibitem{Moser} J.\ Moser, \emph{New aspects in the theory of stability of Hamiltonian systems,} 
Comm. Pure Appl. Math. 11 (1958) 81-114. 

\bibitem{langmuir} I.\,Langmuir, \emph{The structure of the Helium Atom}, Phys.\,Rev.
\textbf{17} (1921), 339--353. 
\bibitem{levi-civita}T.\,Levi-Civita, \emph{Sur la r\'egularisation du probl\`eme des trois corps},
 Acta Math. \textbf{42} (1920), 99--144.
\bibitem{tanner-richter-rost} G.\,Tanner, K.\,Richter, J.\,Rost, \emph{The
theory of two-electron atoms: Between ground state and complete fragmentation,}
Review of Modern Physics \textbf{72}(2) (2000), 497--544. 
\bibitem{wintgen-richter-tanner} D.\,Wintgen, K.\,Richter, G.\,Tanner,\emph{The Semi-Classical Helium Atom}, in Proceedings of the International School of Physics ``Enrico
Fermi", Course CXIX (1993), 113--143.
\bibitem{zhao} L.\,Zhao, \emph{Shooting for Collinear Periodic Orbits in the
Helium Model}, Preprint. 
\bibitem{zhou} B.\,Zhou, \emph{Iteration formulae for brake orbit and index inequalities
for real pseudoholomorphic curves}, J.\,Fixed Point Theory Appl., 
https://doi.org/10.1007/s11784-021-00928-3 (2022).
 
\end{thebibliography}
\end{document}